\documentclass[a4paper,11pt]{article}
\usepackage{indentfirst,amsmath,latexsym,bm,amssymb,epsfig,multicol}
\usepackage{fullpage}
\usepackage{longtable}
\usepackage{color}
\usepackage{amsfonts}
\usepackage{titlesec}
\usepackage{cite}
\usepackage{bm}
\usepackage{lscape}
\usepackage{mathrsfs}
\newenvironment{proof}[1][Proof]{\textbf{#1.} }
{\ \rule{0.75em}{0.75em}\smallskip}
\newcommand\blfootnote[1]{%
\begingroup
\renewcommand\thefootnote{}\footnote{#1}%
\addtocounter{footnote}{-1}%
\endgroup
}

\begin{document}

\begin{center}
\Large An efficient numerical method for a time-fractional
diffusion equation
\end{center}
\begin{center}
Zhongdi  Cen, Jian Huang$^{*}$\blfootnote{*Corresponding author.
Email: sword@zwu.edu.cn (Jian Huang).} , \ Anbo Le, \  Aimin Xu

Institute of Mathematics, Zhejiang Wanli University, Ningbo,
 China
\end{center}

\begin{quote}\small {\bf Abstract:}
A reaction-diffusion problem with a Caputo time derivative is
considered. An integral discretization scheme on a graded mesh
along with a decomposition of the exact solution is proposed. The
truncation error estimate of the discretization scheme is derived
by using the remainder formula of the linear interpolation and
some inequality estimate techniques. It is proved that the scheme
is second-order convergent by applying a difference analogue of
Gronwall's inequality, which exhibits an enhancement in the
convergence rate compared with the $L1$ schemes. Numerical
experiments are presented to support the theoretical result.



\end{quote}
\setcounter{equation}{0}
\section{Introduction}
This article is prompted by recent
publications~\cite{jgos,Mseojg,msjg} where the authors consider
the following initial-boundary value problem
\begin{eqnarray}
& & D_{t}^{\alpha}u(x,t)+L u(x,t)=f(x,t),\ \ \ \ \ \ \ \ \  (x,t)\in Q:=(0,l)\times (0,T],\\
& & u(x,0)=\phi(x),\ \ \ \ \ \ \ \ \ \ \ \ \ \ \ \ \ \ \ \ \ \, \,
\, \ \ \ \ \ \ \  x\in
[0,l],\\
& & u(0,t)=u(l,t)=0,\ \ \ \ \ \ \ \ \ \ \ \ \ \ \ \ \ \ \ \ \ \  \
t\in(0,T].
\end{eqnarray}
Here $D_{t}^{\alpha}$ denotes a Caputo fractional derivative with
$0<\alpha<1$,
\begin{eqnarray*}
L u(x,t):=-p\frac{\partial^{2}u}{\partial x^{2}}(x,t)+c(x)u(x,t),
\end{eqnarray*}
$p$ is a positive constant, $c\in C[0,l]$ with $c\geq 0$, $f\in
C\left(\bar{Q}\right)$ and $\phi\in C[0,l]$. It is proved
in~\cite{jgos,Mseojg} that under reasonable hypotheses on its
data, problem (1.1)-(1.3) has a unique solution $u$ which
typically exhibits a weak singularity at $t=0$.

In~\cite{Mseojg} a finite difference scheme is proposed, which is
a combination of the standard $L1$ approximation for
$D_{t}^{\alpha}u$ on a graded temporal mesh and the central
difference approximation for $L u$ on a uniform spatial mesh. It
is proved that the scheme converges with order
$O\left(M^{-2}+N^{-\min\left\{2-\alpha,r\alpha\right\}}\right)$,
where $M$ and $N$ are the spatial and temporal discretization
parameters and $r\geq 1$ is the mesh grading. In~\cite{jgos,msjg}
a fitted difference scheme and a preprocessed $L1$ scheme are used
to yield an enhanced convergence rate
$O\left(M^{-2}+N^{-\min\left\{2-\alpha,2r\alpha\right\}}\right)$,
respectively.

In the present paper we construct and analyze an integral
discretization scheme on a graded mesh along with a decomposition
of the exact solution of problem (1.1)-(1.3). The truncation error
estimate of the discretization scheme is derived by using the
remainder formula of the linear interpolation and some inequality
estimate techniques. It is shown that the convergence order of our
scheme is $O\left(M^{-2}+N^{-2}\right)$ by applying a difference
analogue of Gronwall's inequality, which improves the convergence
orders given in~\cite{jgos,Mseojg,msjg}. Numerical experiments are
provided to validate the theoretical result.

{\bf Notation.} Throughout the paper, $C$ will denote a generic
positive constant that is independent of the mesh. Note that $C$
can take different values in different places. We always use the
(pointwise) maximum norm $\left\|\cdot\right\|_{\bar{\Omega}}$,
where $\bar{\Omega}$ is a closed and bounded set.

 \setcounter{equation}{0}
\section{The continuous problem}
As in~\cite[Lemma 1]{jgos}, it is assumed that $\phi\in
C^{4}[0,l]$, $0=\phi(0)=\phi''(0)=\phi(l)=\phi''(l)=f(0,t)=f(l,t)$
for $0\leq t\leq T$, $c\in C^{2}[0,l]$ and $f, f_{x}, f_{xx} \in
C\left(\bar{Q}\right)$, and it is shown that the exact solution
$u$ of problem (1.1)-(1.3) can be decomposed as
\begin{eqnarray}
u(x,t)=z(x)t^{\alpha}+\phi(x)+v(x,t),\ \ \ \ \ \ \ \ \ \ \ \ \ \ \
\ \  \  \left(x,t\right)\in \bar{Q},
\end{eqnarray}
where
\begin{eqnarray}
z(x)=\frac{1}{\Gamma\left(\alpha+1\right)}\left(f(x,0)+p\phi''(x)-c(x)\phi(x)\right),
\end{eqnarray}
and $v(x,t)$ is the solution of the following initial-boundary
value problem
\begin{eqnarray}
& & \left(D_{t}^{\alpha}+L\right)v(x,t)=f(x,t)+g(x,t), \ \ \ \ \ \
\ \
\ \ \ \  \left(x,t\right)\in Q,\ \ \ \ \\
& & v(x,0)=0, \ \ \ \ \ \ \ \ \ \ \ \ \ \ \ \ \ \ \ \ \ \ \ \ \ \
\ \ \ \ \ \ \
 \ \ \ \ \ \ \ \ \ \ x\in [0,l],\\
& & v(0,t)=v(l,t)=0,\ \ \ \ \ \ \ \ \ \ \ \ \ \ \  \ \ \ \ \ \ \ \
\ \ \ \ \  \ \ \, \ \ t\in(0,T],
\end{eqnarray}
where
\begin{eqnarray*}
g(x,t)=-f(x,0)+p z''(x)t^{\alpha}-c(x)z(x)t^{\alpha}.
\end{eqnarray*}
 It is proved
in~\cite[Theorem 1]{jgos}, under extra regularity assumptions
$f_{ttt}\left(\cdot,t\right)\in D\left(L^{1/2}\right)$ and
$\left\|f_{tt}\left(\cdot,t\right)\right\|_{L^{1/2}}+t^{\tilde{\rho}}\left\|f_{ttt}\left(\cdot,t\right)\right\|_{L^{1/2}}\leq
C_{1}$ for all $t\in (0,T]$, where $0<\rho <1$ and $C_{1}$ is a
constant independent of $t$, that $v(x,t)$ satisfies
\begin{eqnarray}
\left|\frac{\partial^{k+\ell}v}{\partial x^{k}\partial
t^{\ell}}\right|\leq C\left(1+t^{2\alpha-\ell}\right), \ \ \ \ \ \
\ \ \ \ \ \ \ \ \ 0\leq k+\ell\leq 4,\ \  \ 0\leq\ell\leq 2
\end{eqnarray}
for all $\left(x,t\right)\in [0,l]\times (0,T]$ and some constant
$C$. The similar bounds have been given in~\cite[Theorem
2.1]{Mseojg}, but with $2\alpha$ replaced by $\alpha$. From these
bounds we know that $v$ is smoother than $u$.

It is shown in~\cite[Lemma 6.2]{kdie}  that the problem
(2.3)-(2.5) can be written as the following equivalent
integral-differential equation with a weakly singular kernel
\begin{eqnarray}
& &
v(x,t)=v(x,0)+\frac{1}{\Gamma(\alpha)}\int_{0}^{t}\left(t-s\right)^{\alpha-1}\left[f(x,s)-L
v(x,s)\right]{\rm d}s+G(x,t),\ \
\ \   (x,t)\in Q,\ \ \ \ \ \  \\
& & v(x,0)=0, \ \ \ \ \ \ \ \ \ \ \ \ \ \ \,
 \ \ \ \ \ \ \ \ \ \ x\in [0,l],\\
& & v(0,t)=v(l,t)=0,\ \ \ \ \ \ \ \ \ \ \ \ \ \  t\in(0,T],
\end{eqnarray}
where
\begin{eqnarray*}
G(x,t)=-\frac{t^{\alpha}}{\Gamma(\alpha+1)}f(x,0)+\frac{t^{2\alpha}}{\Gamma(\alpha)}\left(p
z''(x)-c(x)z(x)\right)B(\alpha+1,\alpha).
\end{eqnarray*}
 In the following we will
discrete this integral-differential equation instead of the
differential equation (1.1)-(1.3).

 \setcounter{equation}{0}
\section{Discretization}
In this section we describe a numerical scheme for the
integral-differential equation (2.7)-(2.9). The numerical scheme
is based on a quadrature rule for the integral term and a central
difference method for the temporal discretization.

Based on the properties of the exact solution $v(x,t)$ we
construct a graded mesh $\Omega^{M,N}:=\Omega^{M}\times
\Omega^{N}$, where $\Omega^{M}=\left\{x_{i}=ih\left| 0\leq i\leq
M, h=i/M\right.\right\}$ and $\Omega^{N}=\left\{t_{j}\left| 0\leq
j\leq N, \triangle t_{j}=t_{j}-t_{j-1}\right.\right\}$ with
\begin{eqnarray}
t_{j}=\left\{
\begin{array}{ll}
T\left(\frac{1}{N}\right)^{2/\alpha}, & j=1,\\
T\left(\frac{1}{N}\right)^{2/\alpha}+T\left(\frac{1}{N}\right)^{3/(2\alpha)}, & j=2,\\
T\left(\frac{1}{N}\right)^{2/\alpha}+T\left(\frac{1}{N}\right)^{3/(2\alpha)}+T\left[1-\left(\frac{1}{N}\right)^{2/\alpha}-
\left(\frac{1}{N}\right)^{3/(2\alpha)}\right]\left(\frac{j-2}{N-2}\right)^{1/\alpha},
& 3\leq j\leq N.
\end{array}
\right.
\end{eqnarray}
On this mesh our discrete scheme is second-order convergent.
Furthermore, this mesh avoids too many mesh points concentrating
around $t=0$ compared with the standard graded mesh
$t_{j}=T\left(\frac{j}{N}\right)^{2/\alpha}$ for $0\leq j\leq N$
as that in~\cite{apet1,apet4,nkms,Mseojg}, which improves the
accuracy.

An approximation to the integral can be obtained by the following
quadrature formula
\begin{eqnarray*}
& & \int_{0}^{t_{j}}\left(t_{j}-s\right)^{\alpha-1}\left[f(x,s)-L
v\left(x,s\right)\right]{\rm
d}s\nonumber\\
& & \approx
\sum_{k=1}^{j}\int_{t_{k-1}}^{t_{k}}\left(t_{j}-s\right)^{\alpha-1}\left[\frac{t_{k}-s}{\triangle
t_{k}}\left(f(x,t_{k-1})-L v(x,t_{k-1})\right)
+\frac{s-t_{k-1}}{\triangle t_{k}}\left(f(x,t_{k})-L
v(x,t_{k})\right)\right]{\rm d}s.
\end{eqnarray*}
Then, we have the following discretization scheme for problem
(2.7)-(2.9):
\begin{eqnarray}
\left\{
\begin{array}{ll}
V^{0}_{i}=0, \\
\left\{
\begin{array}{ll}
\displaystyle{V^{j}_{i}=V^{0}_{i}+\frac{1}{\Gamma(\alpha+1)}\sum_{k=1}^{j}\left\{\triangle
t_{k}\left(t_{j}-t_{k-1}\right)^{\alpha}-\frac{1}{\alpha+1}
\left[\left(t_{j}-t_{k-1}\right)^{\alpha+1} \right.\right.}\\
\displaystyle{\left.\left.
-\left(t_{j}-t_{k}\right)^{\alpha+1}\right]\right\}\frac{f^{k-1}_{i}-L^{M}
V^{k-1}_{i}}{\triangle
t_{k}}+\frac{1}{\Gamma(\alpha+1)}\sum_{k=1}^{j} \left\{-\triangle
t_{k}\left(t_{j}-t_{k}\right)^{\alpha}\right.}\\
\displaystyle{\left. +\frac{1}{\alpha+1}
\left[\left(t_{j}-t_{k-1}\right)^{\alpha+1}-\left(t_{j}-t_{k}\right)^{\alpha+1}\right]\right\}\frac{f^{k}_{i}-L^{M}
V^{k}_{i}}{\triangle
t_{k}}}+G_{i}^{j},\ \ \ \ \ \ \ \ \ 1\leq i<M,\\
 V^{j}_{0}=V^{j}_{M}=0,
\end{array}
\right.\\
{\rm for}\ \ j=1,2,\cdots,N,
\end{array}
\right.
\end{eqnarray}
where $V_{i}^{j}$ is the discrete approximation to the exact
solution $v$ of (2.7)-(2.9) at the mesh point $(x_{i},t_{j})$ and
the discrete operator $L^{M}$ is defined as
\begin{eqnarray}
 L^{M}V_{i}^{j}\equiv
-p\frac{V_{i+1}^{j}-2V_{i}^{j}+V_{i-1}^{j}}{h^{2}}+c_{i}V_{i}^{j}.
\end{eqnarray}

 \setcounter{equation}{0}
\section{Convergence analysis}
Let $w^{j}_{i}=V^{j}_{i}-v(x_{i},t_{j})$, where $V_{i}^{j}$ is the
solution of problem (3.2) and $v(x_{i},t_{j})$ is the solution of
problem (2.7)-(2.9) at the mesh point $\left(x_{i},t_{j}\right)$.
Then, the error $w_{i}^{j}$ satisfies the following equation
\begin{eqnarray}
& &
w_{i}^{j}+\frac{1}{\Gamma(\alpha+1)}\sum_{k=1}^{j}\left\{\left(t_{j}-t_{k-1}\right)^{\alpha}-\frac{1}{\left(\alpha+1\right)\triangle
t_{k}} \left[\left(t_{j}-t_{k-1}\right)^{\alpha+1}
-\left(t_{j}-t_{k}\right)^{\alpha+1}\right]\right\}L^{M}
w^{k-1}_{i}\nonumber\\
& &
+\frac{1}{\Gamma(\alpha+1)}\sum_{k=1}^{j}\left\{-\left(t_{j}-t_{k}\right)^{\alpha}
 +\frac{1}{\left(\alpha+1\right)\triangle
t_{k}}
\left[\left(t_{j}-t_{k-1}\right)^{\alpha+1}-\left(t_{j}-t_{k}\right)^{\alpha+1}\right]\right\}L^{M}
w^{k}_{i}\nonumber\\
& & =R^{j}_{i},\ \ \ \ \ \ \ \ \ \ \ \ \ \ \ \ \ \ \ \ \ \ \ \ \ \ \ \ \ \ \ \ \ \ \ \ \ \ \ \ \ \ 1\leq i<M,\ 1\leq j\leq N,\\
& & w^{0}_{i}=0,\ \ \ \ \ \ \ \ \ \ \ \ \ \ \ \ \ \ \ \ \ \ \ \ \ \ \ \ \ \ \ \ \ \ \  \ \ \ \, \ 1\leq i<M,\\
& & w_{0}^{j}=w_{M}^{j}=0,\ \ \ \ \ \ \ \ \ \ \ \ \ \ \ \ \ \ \ \
\ \ \ \ \ \ \ \ \ \ \ 1\leq j\leq N,
\end{eqnarray}
where
\begin{eqnarray}
R^{j}_{i} & = &
\frac{1}{\Gamma(\alpha)}\sum_{k=1}^{j}\int_{t_{k-1}}^{t_{k}}\left(t_{j}-s\right)^{\alpha-1}\left[\frac{t_{k}-s}{\triangle
t_{k}}f(x_{i},t_{k-1})+\frac{s-t_{k-1}}{\triangle
t_{k}}f(x_{i},t_{k})-f(x_{i},s)\right]{\rm d}s\nonumber\\
& &
+\frac{1}{\Gamma(\alpha)}\sum_{k=1}^{j}\int_{t_{k-1}}^{t_{k}}\left(t_{j}-s\right)^{\alpha-1}L
v(x_{i},s){\rm
d}s\nonumber\\
& &
-\frac{1}{\Gamma(\alpha)}\sum_{k=1}^{j}\int_{t_{k-1}}^{t_{k}}\left(t_{j}-s\right)^{\alpha-1}\left[\frac{t_{k}-s}{\triangle
t_{k}}L v(x_{i},t_{k-1})+\frac{s-t_{k-1}}{\triangle
t_{k}}L v(x_{i},t_{k})\right]{\rm d}s   \nonumber\\
& &
+\frac{1}{\Gamma(\alpha)}\sum_{k=1}^{j}\int_{t_{k-1}}^{t_{k}}\left(t_{j}-s\right)^{\alpha-1}\frac{t_{k}-s}{\triangle
t_{k}}\left[L v(x_{i},t_{k-1})-L^{M}v(x_{i},t_{k-1})\right]{\rm d}s \nonumber\\
& &
+\frac{1}{\Gamma(\alpha)}\sum_{k=1}^{j}\int_{t_{k-1}}^{t_{k}}\left(t_{j}-s\right)^{\alpha-1}\frac{s-t_{k-1}}{\triangle
t_{k}}\left[L v(x_{i},t_{k})-L^{M}v(x_{i},t_{k})\right]{\rm d}s.
\end{eqnarray}

For estimating the truncation error we need the following
remainder formula of Newton interpolation.

{\bf Lemma 4.1} (See~\cite{lhxw})  Assume that
$s_{0},s_{1},\dots,s_{k}\in [a,b]$ are distinct. If $u^{(k)}(s)$
is continuous on $[a,b]$, then
\begin{eqnarray*}
u\left[s_{0},s_{1},\dots,s_{k}\right]& =& \int_{0}^{1}{\rm
d}y_{1}\int_{0}^{y_{1}}{\rm d}y_{2}\dots
\int_{0}^{y_{k-1}}u^{(k)}\left((1-y_{1})s_{0}+(y_{1}-y_{2})s_{1}+\right. \nonumber\\
& & \left. \dots+(y_{k-1}-y_{k})s_{k-1}+y_{k}s_{k}\right){\rm
d}y_{k}.
\end{eqnarray*}

Next we give the following technical results under the graded mesh
$\Omega^{N}$.

{\bf Lemma 4.2} Under some regularity conditions on the data,
there exists a positive constant $C$ independent of $N$ such that
\begin{eqnarray*}
\left|\int_{t_{k-1}}^{t_{k}}\left(t_{j}-s\right)^{\alpha-1}\left[L
v(x,s)-\left(\frac{t_{k}-s}{\triangle t_{k}}L
v(x,t_{k-1})+\frac{s-t_{k-1}}{\triangle t_{k}}L
v(x,t_{k})\right)\right]{\rm d}s\right|\leq CN^{-2}
\end{eqnarray*}
for $k=2,3$ and $j\geq k$.

\begin{proof}
By using the remainder formula of Newton interpolation we have
\begin{eqnarray}
& &
\left|\int_{t_{k-1}}^{t_{k}}\left(t_{j}-s\right)^{\alpha-1}\left[L
v(x,s)-\left(\frac{t_{k}-s}{\triangle t_{k}}L
v(x,t_{k-1})+\frac{s-t_{k-1}}{\triangle t_{k}}L
v(x,t_{k})\right)\right]{\rm d}s\right|\nonumber\\
& & \leq
\int_{t_{k-1}}^{t_{k}}\left(t_{j}-s\right)^{\alpha-1}\left|Lv[x;s,t_{k-1},t_{k}]\left(s-t_{k-1}\right)\left(t_{k}-s\right)\right|{\rm
d}s\nonumber\\
& & \leq
\int_{t_{k-1}}^{t_{k}}\int_{0}^{1}\int_{0}^{y_{1}}\left(t_{j}-s\right)^{\alpha-1}\left(s-t_{k-1}\right)\left(t_{k}-s\right)\nonumber\\
& & \ \ \ \cdot\left|L\frac{\partial^{2}v}{\partial
t^{2}}\left(x,(1-y_{1})s+(y_{1}-y_{2})t_{k-1}+y_{2}t_{k}\right)\right|{\rm
d}y_{2}{\rm d}y_{1}{\rm
d}s\nonumber\\
& & \leq
C\int_{t_{k-1}}^{t_{k}}\int_{0}^{1}\int_{0}^{y_{1}}\left(t_{j}-s\right)^{\alpha-1}\left(s-t_{k-1}\right)\left(t_{k}-s\right)\nonumber\\
& & \ \ \
\cdot\left\{1+\left[(1-y_{1})s+(y_{1}-y_{2})t_{k-1}+y_{2}t_{k}\right]^{2\alpha-2}\right\}{\rm
d}y_{2}{\rm d}y_{1}{\rm d}s
\end{eqnarray}
for $k=2,3$, where we have used (2.6). For $\alpha=\frac{1}{2}$,
from (4.5) we have
\begin{eqnarray}
& & \left|\int_{t_{k-1}}^{t_{k}}\left(t_{j}-s\right)^{-1/2}\left[L
v(x,s)-\left(\frac{t_{k}-s}{\triangle t_{k}}L
v(x,t_{k-1})+\frac{s-t_{k-1}}{\triangle t_{k}}L
v(x,t_{k})\right)\right]{\rm d}s\right|\nonumber\\
& & \leq C\left(\triangle
t_{k}\right)^{2}\int_{t_{k-1}}^{t_{k}}\left(t_{j}-s\right)^{-1/2}{\rm
d}s\int_{0}^{1}{\rm
d}y_{1}\int_{0}^{y_{1}}\left[(1-y_{1})s+(y_{1}-y_{2})t_{k-1}+y_{2}t_{k}\right]^{-1}{\rm
d}y_{2}\nonumber\\
& & \leq C\left(\triangle
t_{k}\right)^{2}\int_{t_{k-1}}^{t_{k}}\left(t_{j}-s\right)^{-1/2}{\rm
d}s\int_{0}^{1}\frac{y_{1}}{\left(1-y_{1}\right)s+y_{1}t_{k-1}}{\rm
d}y_{1}\nonumber\\
& & \leq C\left(\triangle
t_{k}\right)^{2}t_{k-1}^{-1}\int_{t_{k-1}}^{t_{k}}\left(t_{j}-s\right)^{-1/2}{\rm
d}s\nonumber\\
& & \leq C\left(\triangle
t_{k}\right)^{2}t_{k-1}^{-1}\left[\left(t_{j}-t_{k-1}\right)^{1/2}-\left(t_{j}-t_{k}\right)^{1/2}\right]\nonumber\\
& & \leq CN^{-2}
\end{eqnarray}
with $k=2,3$. For $0<\alpha<\frac{1}{2}$ and
$\frac{1}{2}<\alpha<1$, from (4.5) we have
\begin{eqnarray}
& &
\left|\int_{t_{k-1}}^{t_{k}}\left(t_{j}-s\right)^{\alpha-1}\left[L
v(x,s)-\left(\frac{t_{k}-s}{\triangle t_{k}}L
v(x,t_{k-1})+\frac{s-t_{k-1}}{\triangle t_{k}}L
v(x,t_{k})\right)\right]{\rm d}s\right|\nonumber\\
& & \leq \frac{C}{\triangle
t_{k}}\int_{t_{k-1}}^{t_{k}}\int_{0}^{1}\left(t_{j}-s\right)^{\alpha-1}\left(s-t_{k-1}\right)\left(t_{k}-s\right)\nonumber\\
& & \ \ \
\cdot\frac{1}{2\alpha-1}\left\{\left[s+(t_{k}-s)y_{1}\right]^{2\alpha-1}-\left[s-(s-t_{k-1})y_{1}\right]^{2\alpha-1}\right\}{\rm
d}y_{1}{\rm
d}s\nonumber\\
& & \leq
C\left|t_{k}^{2\alpha}-t_{k-1}^{2\alpha}\right|\int_{t_{k-1}}^{t_{k}}\left(t_{j}-s\right)^{\alpha-1}{\rm
d}s\nonumber\\
& & =  C\left|t_{k}^{2\alpha}-t_{k-1}^{2\alpha}\right|\left[\left(t_{j}-t_{k-1}\right)^{\alpha}-\left(t_{j}-t_{k}\right)^{\alpha}\right]\nonumber\\
& & \leq CN^{-2}
\end{eqnarray}
with $k=2,3$, where we have used $n\leq 2(n-1)$ for $n\geq 2$.
Combining (4.6) with (4.7) to complete the proof.
\end{proof}

{\bf Lemma 4.3} There exists a positive constant $C$ independent
of $N$ such that
\begin{eqnarray*}
\sum_{k=4}^{j}\left[\left(t_{j}-t_{k-1}\right)^{\alpha}-\left(t_{j}-t_{k}\right)^{\alpha}\right]\left(\triangle
t_{k}\right)^{2}t_{k-1}^{2\alpha-2}\leq CN^{-2}, \ \ \ \ \ \ \ \ \
\  \  4\leq j\leq N.
\end{eqnarray*}

\begin{proof}
Let $\lceil s\rceil$ denote the smallest positive integer that is
greater than or equal to $s$ for any $s\in \mathbb{R}^{+}$. Then
we have
\begin{eqnarray}
& & \sum_{k=4}^{\lceil
j/2\rceil}\left[\left(t_{j}-t_{k-1}\right)^{\alpha}-\left(t_{j}-t_{k}\right)^{\alpha}\right]\left(\triangle
t_{k}\right)^{2}t_{k-1}^{2\alpha-2} \nonumber\\
& & \leq \sum_{k=4}^{\lceil
j/2\rceil}\alpha\left(t_{j}-t_{k}\right)^{\alpha-1}\left(\triangle t_{k}\right)^{3}t_{k-1}^{2\alpha-2}\nonumber\\
& & \leq
\alpha\left(t_{j}-t_{\lceil j/2\rceil}\right)^{\alpha-1}\sum_{k=4}^{\lceil j/2\rceil}\left(\triangle t_{k}\right)^{3}t_{k-1}^{2\alpha-2}\nonumber\\
& & \leq
C\left(\frac{j-2}{N-2}\right)^{1-1/\alpha}\sum_{k=4}^{\lceil
j/2\rceil}\left[\left(\frac{k-2}{N-2}\right)^{1/\alpha}
-\left(\frac{k-3}{N-2}\right)^{1/\alpha}\right]^{3}\left(\frac{k-3}{N-2}\right)^{2-2/\alpha}\nonumber\\
& & \leq C \left(\frac{j-2}{N-2}\right)^{1-1/\alpha}
\left(\frac{1}{N-2}\right)^{3}\sum_{k=4}^{\lceil
j/2\rceil}\left(\frac{k-2}{N-2}\right)^{3/\alpha-3}
\left(\frac{k-3}{N-2}\right)^{2-2/\alpha}\nonumber\\
& & \leq C\left(\frac{1}{N-2}\right)^{3}\sum_{k=4}^{\lceil
j/2\rceil}\left(\frac{k-2}{j-2}\right)^{1/\alpha-1}\nonumber\\
& & \leq CN^{-2},
\end{eqnarray}
where we have used the mean value theorem and $n\leq 2(n-1)$ for
$n\geq 2$. Moreover, we have
\begin{eqnarray}
& & \sum_{k=\lceil
j/2\rceil+1}^{j}\left[\left(t_{j}-t_{k-1}\right)^{\alpha}-\left(t_{j}-t_{k}\right)^{\alpha}\right]\left(\triangle
t_{k}\right)^{2}t_{k-1}^{2\alpha-2} \nonumber\\
& & \leq \max_{\lceil j/2\rceil+1\leq k\leq j}\left(\triangle
t_{k}\right)^{2}t_{k-1}^{2\alpha-2}\sum_{k=\lceil
j/2\rceil+1}^{j}\left[\left(t_{j}-t_{k-1}\right)^{\alpha}
-\left(t_{j}-t_{k}\right)^{\alpha}\right]\nonumber\\
& & \leq  t_{\lceil j/2\rceil}^{2\alpha-2}\left(t_{j}-t_{\lceil
j/2\rceil}\right)^{\alpha} \max_{\lceil j/2\rceil+1\leq k\leq j}\left(\triangle t_{k}\right)^{2}\nonumber\\
& & \leq C\left(\frac{\lceil
j/2\rceil-2}{N-2}\right)^{2-2/\alpha}\frac{j-2}{N-2} \max_{\lceil
j/2\rceil+1\leq k\leq
j}\left[\left(\frac{k-2}{N-2}\right)^{1/\alpha}
-\left(\frac{k-3}{N-2}\right)^{1/\alpha}\right]^{2}\nonumber\\
& & \leq
C\left(\frac{j-2}{N-2}\right)^{2-2/\alpha}\frac{j-2}{N-2}\left(\frac{j-2}{N-2}\right)^{2/\alpha-2}\left(\frac{1}{N-2}\right)^{2}\nonumber\\
& & \leq CN^{-2},
\end{eqnarray}
where we also have used $n\leq 2(n-1)$ for $n\geq 2$. Combining
(4.8) with (4.9) to complete the proof.
\end{proof}

Now we can give the truncation error estimate of the
discretization scheme.

{\bf Lemma 4.4} Under some regularity conditions on the data,
there exists a positive constant $C$ independently of $M$ and $N$
such that the truncation errors of the discretization scheme (3.2)
satisfy
\begin{eqnarray}
\left|R^{j}_{i}\right|\leq  C\left(M^{-2}+N^{-2}\right),\ \ \ \ \
\ \ \ \ 1\leq i\leq M,\ 1\leq j\leq N.
\end{eqnarray}

\begin{proof} For the analysis of the truncation errors we
distinguish two cases.

{\bf Case I}: $j=1$.

From (4.4) we have
\begin{eqnarray}
\left|R^{1}_{i}\right| & \leq &
\frac{1}{\Gamma(\alpha)}\int_{0}^{t_{1}}\left(t_{1}-s\right)^{\alpha-1}\left|\frac{t_{1}-s}{\triangle
t_{1}}f(x_{i},0)+\frac{s}{\triangle
t_{1}}f(x_{i},t_{1})-f(x_{i},s)\right|{\rm d}s\nonumber\\
& &
+\frac{1}{\Gamma(\alpha)}\int_{0}^{t_{1}}\left(t_{1}-s\right)^{\alpha-1}\left|L
v(x_{i},s)\right|{\rm
d}s\nonumber\\
& &
+\frac{1}{\Gamma(\alpha)}\int_{0}^{t_{1}}\left(t_{1}-s\right)^{\alpha-1}\left[\frac{t_{1}-s}{\triangle
t_{1}}\left|L v(x_{i},0)\right|+\frac{s}{\triangle
t_{1}}\left|L v(x_{i},t_{1})\right|\right]{\rm d}s   \nonumber\\
& &
+\frac{1}{\Gamma(\alpha)}\int_{0}^{t_{1}}\left(t_{1}-s\right)^{\alpha-1}\frac{t_{1}-s}{\triangle
t_{1}}\left|L v(x_{i},0)-L^{M}v(x_{i},0)\right|{\rm d}s \nonumber\\
& &
+\frac{1}{\Gamma(\alpha)}\int_{0}^{t_{1}}\left(t_{1}-s\right)^{\alpha-1}\frac{s}{\triangle
t_{1}}\left|L v(x_{i},t_{1})-L^{M}v(x_{i},t_{1})\right|{\rm
d}s\nonumber\\
& \leq & C\left(t_{1}^{\alpha}+t_{1}^{\alpha}M^{-2}\right)\leq
CN^{-2},
\end{eqnarray}
where we have used the assumptions for $f$, (2.6), (3.1) and a
Taylor expansion for $v(x,\cdot)$ about $x_{i}$. From this we
conclude that the lemma holds true for Case I.

{\bf Case II}: $1<j\leq N$.

We decompose the truncation error into two components as follows
\begin{eqnarray}
R^{j}_{i} = R^{j,1}_{i}+R^{j,2}_{i},
\end{eqnarray}
where
\begin{eqnarray}
R^{j,1}_{i} & = &
\frac{1}{\Gamma(\alpha)}\int_{0}^{t_{1}}\left(t_{j}-s\right)^{\alpha-1}\left[\frac{t_{1}-s}{\triangle
t_{1}}f(x_{i},0)+\frac{s}{\triangle
t_{1}}f(x_{i},t_{1})-f(x_{i},s)\right]{\rm d}s\nonumber\\
& &
+\frac{1}{\Gamma(\alpha)}\int_{0}^{t_{1}}\left(t_{j}-s\right)^{\alpha-1}L
v(x_{i},s){\rm
d}s\nonumber\\
& &
-\frac{1}{\Gamma(\alpha)}\int_{0}^{t_{1}}\left(t_{j}-s\right)^{\alpha-1}\left[\frac{t_{1}-s}{\triangle
t_{1}}L v(x_{i},0)+\frac{s}{\triangle
t_{1}}L v(x_{i},t_{1})\right]{\rm d}s   \nonumber\\
& &
+\frac{1}{\Gamma(\alpha)}\int_{0}^{t_{1}}\left(t_{j}-s\right)^{\alpha-1}\frac{t_{1}-s}{\triangle
t_{1}}\left[L v(x_{i},0)-L^{M}v(x_{i},0)\right]{\rm d}s \nonumber\\
& &
+\frac{1}{\Gamma(\alpha)}\int_{0}^{t_{1}}\left(t_{j}-s\right)^{\alpha-1}\frac{s}{\triangle
t_{1}}\left[L v(x_{i},t_{1})-L^{M}v(x_{i},t_{1})\right]{\rm d}s,
\end{eqnarray}
and
\begin{eqnarray}
R^{j,2}_{i} & = &
\frac{1}{\Gamma(\alpha)}\sum_{k=2}^{j}\int_{t_{k-1}}^{t_{k}}\left(t_{j}-s\right)^{\alpha-1}\left[\frac{t_{k}-s}{\triangle
t_{k}}f(x_{i},t_{k-1})+\frac{s-t_{k-1}}{\triangle
t_{k}}f(x_{i},t_{k})-f(x_{i},s)\right]{\rm d}s\nonumber\\
& &
+\frac{1}{\Gamma(\alpha)}\sum_{k=2}^{j}\int_{t_{k-1}}^{t_{k}}\left(t_{j}-s\right)^{\alpha-1}L
v(x_{i},s){\rm
d}s\nonumber\\
& &
-\frac{1}{\Gamma(\alpha)}\sum_{k=2}^{j}\int_{t_{k-1}}^{t_{k}}\left(t_{j}-s\right)^{\alpha-1}\left[\frac{t_{k}-s}{\triangle
t_{k}}L v(x_{i},t_{k-1})+\frac{s-t_{k-1}}{\triangle
t_{k}}L v(x_{i},t_{k})\right]{\rm d}s   \nonumber\\
& &
+\frac{1}{\Gamma(\alpha)}\sum_{k=2}^{j}\int_{t_{k-1}}^{t_{k}}\left(t_{j}-s\right)^{\alpha-1}\frac{t_{k}-s}{\triangle
t_{k}}\left[L v(x_{i},t_{k-1})-L^{M}v(x_{i},t_{k-1})\right]{\rm d}s \nonumber\\
& &
+\frac{1}{\Gamma(\alpha)}\sum_{k=2}^{j}\int_{t_{k-1}}^{t_{k}}\left(t_{j}-s\right)^{\alpha-1}\frac{s-t_{k-1}}{\triangle
t_{k}}\left[L v(x_{i},t_{k})-L^{M}v(x_{i},t_{k})\right]{\rm d}s.
\end{eqnarray}
Similarly to Case I, from (4.13) we have
\begin{eqnarray}
\left|R^{j,1}_{i}\right| & \leq &
\frac{1}{\Gamma(\alpha)}\int_{0}^{t_{1}}\left(t_{j}-s\right)^{\alpha-1}\left|\frac{t_{1}-s}{\triangle
t_{1}}f(x_{i},0)+\frac{s}{\triangle
t_{1}}f(x_{i},t_{1})-f(x_{i},s)\right|{\rm d}s\nonumber\\
& &
+\frac{1}{\Gamma(\alpha)}\int_{0}^{t_{1}}\left(t_{j}-s\right)^{\alpha-1}\left|L
v(x_{i},s)\right|{\rm
d}s\nonumber\\
& &
-\frac{1}{\Gamma(\alpha)}\int_{0}^{t_{1}}\left(t_{j}-s\right)^{\alpha-1}\left|\frac{t_{1}-s}{\triangle
t_{1}}L v(x_{i},0)+\frac{s}{\triangle
t_{1}}L v(x_{i},t_{1})\right|{\rm d}s   \nonumber\\
& &
+\frac{1}{\Gamma(\alpha)}\int_{0}^{t_{1}}\left(t_{j}-s\right)^{\alpha-1}\frac{t_{1}-s}{\triangle
t_{1}}\left|L v(x_{i},0)-L^{M}v(x_{i},0)\right|{\rm d}s \nonumber\\
& &
+\frac{1}{\Gamma(\alpha)}\int_{0}^{t_{1}}\left(t_{j}-s\right)^{\alpha-1}\frac{s}{\triangle
t_{1}}\left|L v(x_{i},t_{1})-L^{M}v(x_{i},t_{1})\right|{\rm
d}s\nonumber\\
& \leq &
C\left[t_{j}^{\alpha}-\left(t_{j}-t_{1}\right)^{\alpha}\right]\left[\left(\triangle
t_{1}\right)^{2}+1+M^{-2}\right]\nonumber\\
& \leq & Ct_{1}\left(t_{j}-t_{1}\right)^{\alpha-1}
 \leq  Ct_{1}^{\alpha}\leq CN^{-2},
\end{eqnarray}
where we also have used the assumptions for $f$, (2.6), (3.1), the
remainder formula of the linear interpolation for $f(\cdot,t)$ and
a Taylor expansion for $v(x,\cdot)$ about $x_{i}$. From (4.14) we
have
\begin{eqnarray}
\left|R^{j,2}_{i}\right| & \leq &
\frac{1}{\Gamma(\alpha)}\sum_{k=2}^{j}\int_{t_{k-1}}^{t_{k}}\left(t_{j}-s\right)^{\alpha-1}\left|\frac{t_{k}-s}{\triangle
t_{k}}f(x_{i},t_{k-1})+\frac{s-t_{k-1}}{\triangle
t_{k}}f(x_{i},t_{k})-f(x_{i},s)\right|{\rm d}s\nonumber\\
& &
+\frac{1}{\Gamma(\alpha)}\sum_{k=2}^{j}\int_{t_{k-1}}^{t_{k}}\left(t_{j}-s\right)^{\alpha-1}\left|L
v(x_{i},s)-\frac{t_{k}-s}{\triangle t_{k}}L
v(x_{i},t_{k-1})-\frac{s-t_{k-1}}{\triangle t_{k}}L
v(x_{i},t_{k})\right|{\rm
d}s  \nonumber\\
& &
+\frac{1}{\Gamma(\alpha)}\sum_{k=2}^{j}\int_{t_{k-1}}^{t_{k}}\left(t_{j}-s\right)^{\alpha-1}\frac{t_{k}-s}{\triangle
t_{k}}\left|L v(x_{i},t_{k-1})-L^{M}v(x_{i},t_{k-1})\right|{\rm d}s \nonumber\\
& &
+\frac{1}{\Gamma(\alpha)}\sum_{k=2}^{j}\int_{t_{k-1}}^{t_{k}}\left(t_{j}-s\right)^{\alpha-1}\frac{s-t_{k-1}}{\triangle
t_{k}}\left|L v(x_{i},t_{k})-L^{M}v(x_{i},t_{k})\right|{\rm
d}s\nonumber\\
& \leq &
\frac{1}{\Gamma(\alpha)}\sum_{k=2}^{j}\left|\frac{\partial^{2}f}{\partial
t^{2}}(x_{i},\xi_{k})\right|\left(\triangle
t_{k}\right)^{2}\int_{t_{k-1}}^{t_{k}}\left(t_{j}-s\right)^{\alpha-1}{\rm
d}s\nonumber\\
& &
+\frac{1}{\Gamma(\alpha)}\sum_{k=4}^{j}\int_{t_{k-1}}^{t_{k}}\left(t_{j}-s\right)^{\alpha-1}\left|L
v(x_{i},s)-\frac{t_{k}-s}{\triangle t_{k}}L
v(x_{i},t_{k-1})-\frac{s-t_{k-1}}{\triangle t_{k}}L
v(x_{i},t_{k})\right|{\rm
d}s  \nonumber\\
& &
+CN^{-2}+CM^{-2}\sum_{k=2}^{j}\int_{t_{k-1}}^{t_{k}}\left(t_{j}-s\right)^{\alpha-1}{\rm
d}s\nonumber\\
& \leq &
C\sum_{k=2}^{j}\left[\left(t_{j}-t_{k-1}\right)^{\alpha}-\left(t_{j}-t_{k}\right)^{\alpha}\right]
\left[\left(\triangle
t_{k}\right)^{2}+M^{-2}\right]+CN^{-2}\nonumber\\
& &
+\frac{1}{\Gamma(\alpha)}\sum_{k=4}^{j}\left|L\frac{\partial^{2}v}{\partial
t^{2}}(x,\eta_{k})\right|\left(\triangle
t_{k}\right)^{2}\int_{t_{k-1}}^{t_{k}}\left(t_{j}-s\right)^{\alpha-1}{\rm
d}s\nonumber\\
& \leq &
C\left(M^{-2}+N^{-2}\right)\left(t_{j}-t_{1}\right)^{\alpha}+CN^{-2}+C\sum_{k=4}^{j}\left[\left(t_{j}-t_{k-1}\right)^{\alpha}-\left(t_{j}-t_{k}\right)^{\alpha}\right]
\left(\triangle
t_{k}\right)^{2}t_{k-1}^{2\alpha-2}\nonumber\\
& \leq &  C\left(M^{-2}+N^{-2}\right),
\end{eqnarray}
where we have used the remainder formula of the linear
interpolation for $v(\cdot,t)$ and $f(\cdot,t)$ with $\xi_{k},
\eta_{k}\in (t_{k-1},t_{k})$, the bounds on $v(x,t)$ and its
derivatives given by (2.6), the assumptions for $f(x,t)$, Lemmas
4.2 and 4.3. Therefore, from (4.12), (4.15) and (4.16) we conclude
that the lemma also holds true for Case II.
\end{proof}

Next we give the error estimates for the discretization scheme.

{\bf Theorem 4.5} Let $v(x,t)$ be the solution of problem
(2.7)-(2.9) and $V$ be the solution of problem (3.2). Then, under
some regularity conditions on the data, we have the following
error estimate
\begin{eqnarray}
\left\|V-v\right\|_{\Omega^{M,N}}\leq C\left(M^{-2}+N^{-2}\right),
\end{eqnarray}
where $C$ is a positive constant independent of $M$ and $N$.

\begin{proof}
From (4.1) we have
\begin{eqnarray}
w^{j}_{i} & = & \left(I+\frac{\left(\triangle
t_{j}\right)^{\alpha}}{\Gamma(\alpha+2)}L^{M}\right)^{-1}R^{j}_{i}-\frac{1}{\Gamma(\alpha+1)}\sum_{k=1}^{j}
\left\{\left(t_{j}-t_{k-1}\right)^{\alpha}\right.\nonumber\\
& & \left. -\frac{1}{(\alpha+1)\triangle t_{k}}\left[\left(t_{j}-t_{k-1}\right)^{\alpha+1}-\left(t_{j}-t_{k}\right)^{\alpha+1}\right]\right\}\left(I+\frac{\left(\triangle
t_{j}\right)^{\alpha}}{\Gamma(\alpha+2)}L^{M}\right)^{-1}L^{M}w^{k-1}_{i}\nonumber\\
& & -\frac{1}{\Gamma(\alpha+1)}\sum_{k=1}^{j-1}
\left\{-\left(t_{j}-t_{k}\right)^{\alpha}
+\frac{1}{(\alpha+1)\triangle
t_{k}}\left[\left(t_{j}-t_{k-1}\right)^{\alpha+1}-\left(t_{j}-t_{k}\right)^{\alpha+1}\right]\right\}\nonumber\\
& & \cdot\left(I+\frac{\left(\triangle
t_{j}\right)^{\alpha}}{\Gamma(\alpha+2)}L^{M}\right)^{-1}L^{M}w^{k}_{i}.
\end{eqnarray}

It is easy to see that the operator $\left(I+\frac{\left(\triangle
t_{j}\right)^{\alpha}}{\Gamma(\alpha+2)}L^{M}\right)$ satisfies a
discrete maximum principle, and consequently
\begin{eqnarray}
\left\|\left(I+\frac{\left(\triangle
t_{j}\right)^{\alpha}}{\Gamma(\alpha+2)}L^{M}\right)^{-1}\right\|_{\Omega^{M}}\leq
1, \ \ \ \ \ \ \ \ \ \ \ \   1\leq j\leq N.
\end{eqnarray}
Furthermore, applying the result proved in Palencia~\cite{cpal} we
have
\begin{eqnarray}
\left\|\left(I+\frac{\left(\triangle
t_{j}\right)^{\alpha}}{\Gamma(\alpha+2)}L^{M}\right)^{-1}L^{M}\right\|_{\Omega^{M}}\leq
d_{j}, \ \ \ \ \ \ \ \ \ \ \ \ \ 1\leq j\leq N,
\end{eqnarray}
since $\left(1+\frac{\left(\triangle
t_{j}\right)^{\alpha}}{\Gamma(\alpha+2)}y\right)^{-1} y$ is a
rational A-acceptable function, where $d_{j}$ is a positive
constant. The analogous problems have been discussed
in~\cite{ccjgj,mktak}.

 Therefore, from (4.18)-(4.20) we can
obtain
\begin{eqnarray}
\left\|w^{j}\right\|_{\Omega^{M}}\leq
z_{j}+d_{j}\sum_{k=1}^{j-1}q_{k}\left\|w^{k}\right\|_{\Omega^{M}},
\end{eqnarray}
where
\begin{eqnarray*}
z_{k}& =&
\left\|R^{k}\right\|_{\Omega^{M}},\\
q_{k} & = & \frac{1}{\Gamma(\alpha+2)}\left\{\frac{1}{\triangle
t_{k}}\left[\left(t_{j}-t_{k-1}\right)^{\alpha+1}-\left(t_{j}-t_{k}\right)^{\alpha+1}\right] \right.\\
& & \left. -\frac{1}{\triangle
t_{k+1}}\left[\left(t_{j}-t_{k}\right)^{\alpha+1}-\left(t_{j}-t_{k+1}\right)^{\alpha+1}\right]\right\},
\ \ \ \ \ \ \  1\leq k\leq j.
\end{eqnarray*}
Then applying the discrete analogue of Gronwall's
inequality~\cite[Theorem 3]{dwjw}, we have
\begin{eqnarray}
\left\|w^{j}\right\|_{\Omega^{M}} \leq
z_{j}+d_{j}\prod_{m=1}^{j-1}\left(1+d_{m}q_{m}\right)\cdot\sum_{k=1}^{j-1}\left[z_{k}q_{k}\prod_{m=1}^{k}\left(1+d_{m}q_{m}\right)^{-1}\right]
\end{eqnarray}
for $1\leq j\leq N$. Lemma 4.4 implies
\begin{eqnarray}
0<z_{k}\leq C\left(M^{-2}+N^{-2}\right),\ \ \ \ \ \ \ \ \ \ \ \  \
1\leq k\leq j.
\end{eqnarray}
Furthermore, we have
\begin{eqnarray}
\sum_{k=1}^{j-1}q_{k} & = &
\frac{1}{\Gamma(\alpha+2)}\sum_{k=1}^{j-1}\left\{\frac{1}{\triangle
t_{k}}\left[\left(t_{j}-t_{k-1}\right)^{\alpha+1}-\left(t_{j}-t_{k}\right)^{\alpha+1}\right] \right.\nonumber\\
& & \left. -\frac{1}{\triangle
t_{k+1}}\left[\left(t_{j}-t_{k}\right)^{\alpha+1}-\left(t_{j}-t_{k+1}\right)^{\alpha+1}\right]\right\}\nonumber\\
& = & \frac{1}{\Gamma(\alpha+1)}\sum_{k=1}^{j-1}\left[\left(t_{j}-\mu_{k}\right)^{\alpha}-\left(t_{j}-\mu_{k+1}\right)^{\alpha}\right]\nonumber\\
& = &
\frac{1}{\Gamma(\alpha+1)}\left[\left(t_{j}-\mu_{1}\right)^{\alpha}-\left(t_{j}-\mu_{j}\right)^{\alpha}\right],
\end{eqnarray}
where we have used the mean value theorem with $\mu_{k}\in
(t_{k-1},t_{k})$. Thus we have
\begin{eqnarray}
& &
\sum_{k=1}^{j-1}\left[z_{k}q_{k}\prod_{m=1}^{k}\left(1+d_{m}q_{m}\right)^{-1}\right]
\leq
CN^{-2}\sum_{k=1}^{j-1}q_{k}\nonumber\\
& & = C\left(M^{-2}+N^{-2}\right)\left[\left(t_{j}-\mu_{1}\right)^{\alpha}-\left(t_{j}-\mu_{j}\right)^{\alpha}\right]\nonumber\\
& & \leq C\left(M^{-2}+N^{-2}\right),
\end{eqnarray}
and
\begin{eqnarray}
\prod_{m=1}^{j-1}\left(1+d_{m}q_{m}\right) \leq
\exp\left(\sum_{m=1}^{j-1}d_{m}q_{m}\right)
 \leq
\exp\left(C\left[\left(t_{j}-\mu_{1}\right)^{\alpha}-\left(t_{j}-\mu_{j}\right)^{\alpha}\right]\right)\leq
C,
\end{eqnarray}
where we have used (4.23)-(4.24) and the inequality
$\left(1+y\right)\leq e^{y}$ for $y\geq -1$. Hence, combining
(4.25)-(4.26) with (4.22) we can obtain
\begin{eqnarray*}
\left\|w^{j}\right\|_{\Omega^{M}} \leq
C\left(M^{-2}+N^{-2}\right), \ \ \ \ \ \ \ \ \ \ \ \ \ \ 1\leq
j\leq N.
\end{eqnarray*}
From this we complete the proof.
\end{proof}

Then, our approximation $U_{i}^{j}$ of $u(x_{i},t_{j})$ can be
obtained from (2.1)
\begin{eqnarray}
U_{i}^{j}=z(x_{i})t_{j}^{\alpha}+\phi(x_{i})+V_{i}^{j},\ \ \ \ \ \
\ \ 0\leq i\leq M,\ 0\leq j\leq N.
\end{eqnarray}
Therefore, from (4.27) and Theorem 4.5 we have
\begin{eqnarray}
\left\|U-u\right\|_{\Omega^{M,N}}\leq C\left(M^{-2}+N^{-2}\right),
\end{eqnarray}
which improves the convergence orders given
in~\cite{jgos,Mseojg,msjg}. There are two reasons for an
enhancement in the convergence rate. The first one is that the
fractional differential equation is transformed into an equivalent
integral-differential equation which reduces the singularity of
the integrand function. The other reason is that the decomposition
of the exact solution is used and the remainder term $v$ is
smoother than $u$.

\setcounter{equation}{0}
\section{Numerical experiments}
In this section we verify experimentally the theoretical results
obtained in the preceding section. Error estimates and convergence
rates for the discrete scheme are presented for the following
example which has been given in~\cite{jgos,msjg}.

{\bf Example} Fractional differential equation with
non-homogeneous boundary conditions:
\begin{eqnarray*}
& & D_{t}^{\alpha}u-\frac{\partial^{2}u}{\partial x^{2}}=f(x,t),\
\ \ \ \ \ \ \ \  \ \ \ \
(x,t)\in (0,\pi)\times (0,1],\\
& & u(x,0)=\sin x,\ \ \ \ \ \ \ \ \
\ \ \ \ \ \ \ \ \ \ \ \ \  x\in (0,\pi),\\
& & u(0,t)=u(1,t)=0,\ \ \ \ \ \ \ \ \ \ \ \ \ \  \ t\in (0,1].
\end{eqnarray*}
The function $f(x,t)$ is chosen such that the exact solution is
$u(x,t)=\left[E_{\alpha}\left(-t^{\alpha}\right)+t^{3}\right]\sin
x$, where $E_{\alpha}\left(\cdot\right)$ is the classical
Mittag-Leffler function. The solution $u(x,t)$ has a typical weak
singularity at $t=0$ (see~\cite{jgos,msjg}).

The maximum error is denoted by
\begin{eqnarray*}
e^{M,N}=\left\|U-u\right\|_{\bar{\Omega}^{M,N}},
\end{eqnarray*}
and the corresponding convergence rate is computed by
\begin{eqnarray*}
rate^{M,N}=\log_{2}\left(\frac{e^{M,N}}{e^{2M,2N}}\right)
\end{eqnarray*}
for the discrete scheme (3.2). The error estimates and convergence
rates in our computed solutions are listed in Table 1. Table 1
shows that the computed solution converges to the exact solution
with second-order accuracy and the numerical results do not depend
strongly on the value of $\alpha$, which supports the convergence
estimate of Theorem 4.5.

For comparison we also use the standard $L1$ scheme~\cite{Mseojg}
with $r=\left(2-\alpha\right)/\alpha$ (optimal choice) and the
preprocessed $L1$ scheme~\cite{msjg} with
$r=\left(2-\alpha\right)/(2\alpha)$ (optimal choice) to compute
this example. The numerical results are presented in Table 2. From
Tables 1 and 2 we confirm that our method proposed in this paper
is more accurate and robust than the $L1$ scheme and the
preprocessed $L1$ scheme.

\begin{table} \caption{Error estimates $e^{M,N}$ and convergence rates $rate^{M,N}$ of the scheme (3.2) for Example}
\begin{center}
\begin{tabular}{c}
 \hline
$M=N$\ \ \ \ \ \ \ \ \  $64$\ \ \ \ \ \ \ \ \ \  \ \ \  \ \ $128$
\ \ \ \ \ \ \ \
\ \ \ \ \ \ \   $256$\ \ \ \ \ \ \ \ \ \ \ \ \  \ \   $512$\ \ \ \ \ \ \ \ \ \ \ \  \ \ \ \  $1024$\ \ \ \ \\
 \hline
$\alpha=0.2$\ \ \ \ \  1.0185e-3 \  \ \ \ \ \   2.7198e-4 \ \ \ \
\ \ \ 7.2032e-5 \ \ \ \ \ \ \   1.8931e-5
\ \ \ \ \ \ \  4.9100e-6 \\
\ \ \ \ \ \ \ \ \ \ \ \ \ \ \ \ \ \ \ \ \ \   1.905  \ \ \ \ \ \ \
\ \ \ \ 1.917 \ \ \ \ \ \ \ \ \ \ \ \  1.928 \ \ \ \ \ \ \ \
\ \ \ \ \  1.947  \ \ \ \ \ \ \ \ \ \ \ \ \ \ \ \  -  \ \ \ \ \ \ \ \ \ \ \  \\
$\alpha=0.4$\ \ \ \ \   4.7052e-4 \  \ \ \ \ \   1.1803e-4 \ \ \ \
\ \ \ 2.9727e-5 \ \ \ \ \ \ \  7.4922e-6
\ \ \ \ \ \ \  1.8869e-6 \\
\ \ \ \ \ \ \ \ \ \ \ \ \ \ \ \ \ \ \ \ \ \   1.995  \ \ \ \ \ \ \
\ \ \ \ 1.989 \ \ \ \ \ \ \ \ \ \ \ \  1.988  \ \ \ \ \ \ \ \
\ \ \ \ \  1.989  \ \ \ \ \ \ \ \ \ \ \ \ \ \ \ \  -  \ \ \ \ \ \ \ \ \ \ \  \\
$\alpha=0.6$\ \ \ \ \  2.7573e-4 \  \ \ \ \ \   6.8004e-5 \ \ \ \
\ \ \ 1.6902e-5 \ \ \ \ \ \ \  4.2153e-6
\ \ \ \ \ \ \  1.0530e-6 \\
\ \ \ \ \ \ \ \ \ \ \ \ \ \ \ \ \ \ \ \ \ \  2.020  \ \ \ \ \ \ \
\ \ \ \ 2.008 \ \ \ \ \ \ \ \ \ \ \ \  2.003   \ \ \ \ \ \ \ \
\ \ \ \ \  2.001  \ \ \ \ \ \ \ \ \ \ \ \ \ \ \ \  -  \ \ \ \ \ \ \ \ \ \ \  \\
$\alpha=0.8$\ \ \ \ \  1.8272e-4  \  \ \ \ \ \  4.4962e-5 \ \ \ \
\ \ \ 1.1153e-5 \ \ \ \ \ \ \  2.7776e-6
\ \ \ \ \ \ \  6.9309e-7  \\
\ \ \ \ \ \ \ \ \ \ \ \ \ \ \ \ \ \ \ \ \ \   2.023   \ \ \ \ \ \
\ \ \ \ \ 2.011 \ \ \ \ \ \ \ \ \ \ \ \  2.006 \ \ \ \ \ \ \ \
\ \ \ \ \  2.003  \ \ \ \ \ \ \ \ \ \ \ \ \ \ \ \  -  \ \ \ \ \ \ \ \ \ \ \  \\
\hline
\end{tabular}
\end{center}
\end{table}
\begin{table} \caption{Error estimates $e^{M,N}$ and convergence rates $rate^{M,N}$ of
the standard $L1$ scheme ($L1$)~\cite{Mseojg} and the preprocessed
$L1$ scheme (P$L1$)~\cite{msjg}
 with optimal $r$ for Example}
\begin{center}
\begin{tabular}{c}
 \hline
$M=N$\ \ \ \ \ \ \ \ \ \ \ \  \ \ \ \ \  \  $64$\ \ \ \ \ \ \ \ \
\ \ \ \ \ \ \ $128$ \ \ \ \ \ \ \ \ \
\ \ \ \ \ \    $256$\ \ \ \ \ \ \ \ \ \ \ \ \  \ \   $512$\ \ \ \ \ \ \ \ \ \ \ \ \   \ \   $1024$ \ \ \ \\
 \hline
$\alpha=0.2$\ \ \ \ \ $L1$ \ \  \  4.5112e-3 \  \ \ \ \ \
1.3940e-3 \ \ \ \ \ \ \ \ 3.6266e-4 \  \ \ \ \ \   2.3831e-4
\ \ \ \ \ \ \  2.6706e-4 \\
\ \ \ \ \ \ \ \ \ \ \ \ \ \ \ \ \ \ \ \  \ \ \ \ \ \ \ \ \ \ \
1.694 \ \ \ \ \ \ \ \ \ \ \ \  1.943 \ \ \ \ \ \ \ \ \ \ \ \ \
0.606 \ \ \ \ \ \
\ \ \ \ \ \  -0.164  \ \ \ \ \ \ \ \ \ \ \ \ \ \ \  -  \ \ \ \ \ \ \ \ \ \ \  \\
\ \ \ \ \ \ \ \ \ \ \ \ \ \  P$L1$ \ \    1.6443e-3 \  \ \ \ \ \
5.2018e-4 \ \ \ \ \ \ \ \ 1.6109e-4  \ \ \ \ \ \   4.9137e-5
\ \ \ \ \ \ \  1.4823e-5 \\
\ \ \ \ \ \ \ \ \ \ \ \ \ \ \ \ \ \ \ \  \ \ \ \ \ \ \ \ \ \ \
 1.660 \ \ \ \ \ \ \ \ \ \ \ \  1.691 \ \ \ \ \ \ \ \ \ \ \ \ \
1.713 \ \ \ \ \ \
\ \ \ \ \ \  1.729   \ \ \ \ \ \ \ \ \ \ \ \ \ \ \  -  \ \ \ \ \ \ \ \ \ \ \  \\
$\alpha=0.4$\ \ \ \ \ $L1$ \ \  \  4.6180e-3 \  \ \ \ \ \
1.6175e-3 \ \ \ \ \ \ \ \ 5.5659e-4 \  \ \ \ \ \   1.8926e-4
\ \ \ \ \ \ \  6.3823e-5 \\
\ \ \ \ \ \ \ \ \ \ \ \ \ \ \ \ \ \ \ \  \ \ \ \ \ \ \ \ \ \ \
1.514 \ \ \ \ \ \ \ \ \ \ \ \  1.539 \ \ \ \ \ \ \ \ \ \ \ \ \
1.556 \ \ \ \ \ \
\ \ \ \ \ \  1.568   \ \ \ \ \ \ \ \ \ \ \ \ \ \ \  -  \ \ \ \ \ \ \ \ \ \ \  \\
\ \ \ \ \ \ \ \ \ \ \ \ \ \  P$L1$ \ \    1.6527e-3 \  \ \ \ \ \
5.5897e-4 \ \ \ \ \ \ \ \ 1.8773e-4  \ \ \ \ \ \   6.2742e-5
\ \ \ \ \ \ \  2.0897e-5 \\
\ \ \ \ \ \ \ \ \ \ \ \ \ \ \ \ \ \ \ \  \ \ \ \ \ \ \ \ \ \ \
1.564 \ \ \ \ \ \ \ \ \ \ \ \  1.574 \ \ \ \ \ \ \ \ \ \ \ \ \
1.581 \ \ \ \ \ \
\ \ \ \ \ \  1.586   \ \ \ \ \ \ \ \ \ \ \ \ \ \ \  -  \ \ \ \ \ \ \ \ \ \ \  \\
$\alpha=0.6$\ \ \ \ \ $L1$ \ \  \  6.2359e-3 \  \ \ \ \ \
2.4091e-3 \ \ \ \ \ \ \ \ 9.2427e-4 \  \ \ \ \ \  3.5303e-4
\ \ \ \ \ \ \  1.3446e-4 \\
\ \ \ \ \ \ \ \ \ \ \ \ \ \ \ \ \ \ \ \  \ \ \ \ \ \ \ \ \ \ \
1.372 \ \ \ \ \ \ \ \ \ \ \ \  1.382 \ \ \ \ \ \ \ \ \ \ \ \ \
1.389 \ \ \ \ \ \
\ \ \ \ \ \  1.393  \ \ \ \ \ \ \ \ \ \ \ \ \ \ \  -  \ \ \ \ \ \ \ \ \ \ \  \\
\ \ \ \ \ \ \ \ \ \ \ \ \ \  P$L1$ \ \    2.5219e-3 \  \ \ \ \ \
9.5577e-4 \ \ \ \ \ \ \ \ 3.6218e-4  \ \ \ \ \ \   1.3723e-4
\ \ \ \ \ \ \  5.1999e-5 \\
\ \ \ \ \ \ \ \ \ \ \ \ \ \ \ \ \ \ \ \  \ \ \ \ \ \ \ \ \ \ \
1.400 \ \ \ \ \ \ \ \ \ \ \ \  1.400 \ \ \ \ \ \ \ \ \ \ \ \ \
1.400 \ \ \ \ \ \
\ \ \ \ \ \   1.400   \ \ \ \ \ \ \ \ \ \ \ \ \ \ \  -  \ \ \ \ \ \ \ \ \ \ \  \\
$\alpha=0.8$\ \ \ \ \ $L1$ \ \  \  1.0663e-2 \  \ \ \ \ \
4.6714e-3 \ \ \ \ \ \ \ \ 2.0426e-3 \  \ \ \ \ \   8.9194e-4
\ \ \ \ \ \ \  3.8915e-4 \\
\ \ \ \ \ \ \ \ \ \ \ \ \ \ \ \ \ \ \ \  \ \ \ \ \ \ \ \ \ \ \
1.191 \ \ \ \ \ \ \ \ \ \ \ \  1.193 \ \ \ \ \ \ \ \ \ \ \ \ \
1.195 \ \ \ \ \ \
\ \ \ \ \ \ 1.197  \ \ \ \ \ \ \ \ \ \ \ \ \ \ \  -  \ \ \ \ \ \ \ \ \ \ \  \\
\ \ \ \ \ \ \ \ \ \ \ \ \ \  P$L1$ \ \    5.9732e-3 \  \ \ \ \ \
2.6142e-3 \ \ \ \ \ \ \ \ 1.1449e-3  \ \ \ \ \ \   5.0145e-4
\ \ \ \ \ \ \  2.1957e-4 \\
\ \ \ \ \ \ \ \ \ \ \ \ \ \ \ \ \ \ \ \  \ \ \ \ \ \ \ \ \ \ \
1.192 \ \ \ \ \ \ \ \ \ \ \ \ 1.191 \ \ \ \ \ \ \ \ \ \ \ \ \
1.191 \ \ \ \ \ \
\ \ \ \ \ \  1.191   \ \ \ \ \ \ \ \ \ \ \ \ \ \ \  -  \ \ \ \ \ \ \ \ \ \ \  \\
\hline
\end{tabular}
\end{center}
\end{table}

{\bf Acknowledgement.} The work was supported by Humanities and
Social Sciences Planning Fund of Ministry of Education of China
(Grant No. 18YJAZH002), Major humanities and Social Sciences
projects in colleges and universities of Zhejiang (Grant No.
2018GH020), Zhejiang Province Natural Science Foundation (Grant
No. Y17D010024).


\begin{thebibliography}{0}
\bibitem{ccjgj}C. Clavero, J.L. Gracia, J.C. Jorge, High-order numerical methods
for one-dimensional parabolic singularly perturbed problems with
regular layers, {\em Numer. Meth. Part. Differ. Equ.}, 21(1)
(2005) 149-169.
\bibitem{kdie}K. Diethlm, {\em The analysis of fractional differential equations}, in:
Lecture Notes in Mathematics, vol. 2004, Springer, Berlin, 2010.
\bibitem{jgos}J.L. Gracia, E. O'Riordan, and M. Stynes, A fitted scheme for a Caputo initial-boundary value problem,
 {\em J. Sci. Comput.}, 76(1) (2018) 583-609.
\bibitem{lhxw}L.C. Hsu, X.H. Wang, Examples and methods in mathematical
analysis, Higher Education Press, 1983, Page 234 (in Chinese).
 \bibitem{mktak}M.K. Kadalbajoo, L.P. Tripathi, and A. Kumar, A cubic B-spline
collocation method for a numerical solution of the generalized
Black-Scholes equation, {\em Math. Comput. Model.}, 55(3-4) (2012)
1483-1505.
 \bibitem{apet4}M. Kolk, A. Pedas, and E. Tamme, Modified spline collocation
for linear fractional differential equations, {\em J. Comput.
Appl. Math.}, 283 (2015) 28-40.
\bibitem{nkms}N. Kopteva, M. Stynes, An efficient collocation method for a Caputo
two-point boundary value problem, {\em BIT Numer. Math.}, 55
(2015) 1105-1123.
\bibitem{cpal}C. Palencia, A stability result for sectorial operators in Banach
spaces, {\em SIAM J. Numer. Anal.}, 30(5) (1993) 1373-1384.
 \bibitem{apet1}A. Pedas, E. Tamme, Piecewise polynomial collocation
for linear boundary value problems of fractional differential
equations, {\em J. Comput. Appl. Math.}, 236 (2012) 3349-3359.
 \bibitem{Mseojg}M. Stynes, E. O'Riordan, and J.L. Gracia, Error analysis of a
finite difference method on graded meshes for a time-fractional
diffusion equation, {\em SIAM J. Numer. Anal.}, 55(2) (2017)
1057-1079.
\bibitem{msjg}M. Stynes, J.L. Gracia, Preprocessing schemes for fractional-derivative problems to improve their convergence rates,
 {\em Appl. Math. Lett.}, 74 (2017) 187-192.
\bibitem{msjlg1}M. Stynes, J.L. Gracia, A finite difference method
for a two-point boundary value problem with a Caputo fractional
derivative, {\em IMA J. Numer. Anal.}, 35 (2015) 689-721.
\bibitem{dwjw}D. Willett, J.S.W. Wong, On the discrete analogues of some
generalizations of Gronwall's inequality, {\em Monatsh. Math.},
69(4) (1965) 362-367.
\end{thebibliography}
\end{document}